\documentclass[11pt]{article}
\usepackage{fontenc,amsmath,mathrsfs,amsfonts,amsthm,graphicx, verbatim,geroeng,amssymb,array}
\begin{document}
\date{}
\title{Generalized Thue-Morse Continued Fractions}
\author{Gerardo Gonz\'alez Robert}
\maketitle
\begin{abstract}
The Thue-Morse sequence is generalized to the $TM_m$ sequences and two equivalent definitions are given. This generalization leads to transcendental numbers and has Queff\'elec's theorem on Thue-Morse continued fractions as a special case. It is also shown that the theorem of Adamczewski and Bugeaud for palindromic continued fractions cannot be applied in this case. 
\end{abstract}
\section{Introduction}
It is a well known fact that almost all real numbers, in the sense of Lebesgue measure, are transcendental. It is also well known that proving the transcendence of a given real number can be quite hard. Liouville numbers might be the reals whose transcendence is the easiest to prove using the approximation condition provided by Liouville's approximation theorem (see \cite{bugo01}, Theorem 1.2, p.3). In simplest terms, this theorem states that if $\alpha$ can be approximated by a sequence of rational numbers $(p_n/q_n)_{n=0}^{\infty}$ such that the error decreases exponentially in terms of $q_n$, then $\alpha$ is transcendental.

We recall that an irrational real number $\theta$ satisfying $\liminf n\|n\theta\|>0$ is \textbf{badly approximable}, where $\|\cdot\|$ is the distance to the closest integer. The set of badly approximable numbers is denoted $\bad$. It is well known that $\theta \in \bad$ if and only if the continued fraction expansion of $\theta$ has bounded partial quotients (see \cite{bugo01}, p. 11, Theorem 1.9 or \cite{schm}, p. 22, Theorem 5F). In particular, the set $\bad$ is uncountable. Moreover, if $a,b$ are distinct positive integers the set of irrationals whose partial quotients take values in $\left\{a,b\right\}$ is uncountable, so there are infinitely many transcendental numbers in this set. Martine Queff\'elec provided in  \cite{queff} a way for constructing transcendental numbers whose partial quotients take only two values.
\begin{teo01} [Queff\'elec]
Let $a,b$ be two distinct positive integers. Then the real number $\alpha$ whose sequence of partial quotients is the Thue-Morse sequence over the alphabet $\left\{a,b\right\}$ is transcendental.
\label{Quef}
\end{teo01}

(See Section 2 for the definition of Thue-Morse sequence). Adamczewski and Bugeaud obtained in \cite{bugam01} Queff\'elec's result as a corollary of the following
\begin{teo01}[Adamczewski-Bugeaud]
Let $\bfa=(a_n)_{n=1}^{\infty}$ be a sequence of positive integers. If the word $\bfa$ is not eventually periodic and is palindromic, then $[0;a_1,a_2,a_3,\ldots]$ is transcendental.
\label{Adabug}
\end{teo01}

(See Section 2 for the definition of palindromic.) In this paper the Thue-Morse sequences are generalized to what we call the $TM_m$ sequences-- sequences taking values in an alphabet in $m$ letters. We show that such sequence give rise to transendental numbers as in Theorem \ref{Quef}; however, as we shall see, we may not use Theorem \ref{Adabug} to show this, as these sequences will not be palindromic except for the case $m=2$. Explicitly, in Section 4 we prove the following

\begin{teo01}
Let $a_1,a_2,\ldots,a_m$ be $m$ pairwise distinct positive integers. Then the real number $\alpha$ whose sequence of partial quotients is the $TM_m$ sequence over the alphabet $\left\{a_1,\ldots,a_m \right\}$ is transcendental.
\label{Mien}
\end{teo01}

\section{$TM_m$ sequences}

We will adopt the notation $\Na:=\left\{ n\in \Za: n> 0\right\}$, $\Na_0=\Na\cup \left\{0\right\}$ and we recall that $[a_0;a_1,a_2,\ldots]$ represents the real number whose simple continued fraction is given precisely by the sequence $(a_n)_{n=0}^{\infty}$. We will work with words, denoted by $\bfw=w_0\ldots w_n$, on numerical alphabets. We will sometimes use commas when referring to specific words, e.g. $n-1, \, m-1$ instead of $n-1 \, m-1$.

There are several ways to define and generalize the Thue-Morse sequence ($TM_2$). We will define $TM_2$ in terms of the binary expansion of a natural number and its generalization will be given in terms of the finite $m$-ary representation of a natural number. We use the notation of \cite{allouche}: for non-negative integers $n,m$, $m\geq 2$, $(n)_m$ is the word formed by the coefficients of the finite $m$-ary representation of $n$; that is, 
\[
n=\sum_{i}^r c_i m^i \; \Rightarrow \; \quad (n)_m:=c_0c_1\ldots c_r. 
\]
For any finite word $\bfw=w_0w_1\ldots w_r$ over the alphabet $\left\{0,1,\ldots,m-1\right\}$, $[\bfw]_m=\sum_{i=0}^r w_im^i$, $|\bfw|=r+1$, and $\bfw_{[j_1:j_2]}=w_{j_1}w_{j_1+1}\ldots w_{j_2-1}$. Also, $\Sigma_m:=\left\{1,2,\ldots,m\right\}$.
\begin{def01}
Let $m\geq 2$ be an integer. Define the sequence $\left(t_{n}{(m)}\right)_{n=0}^{\infty}$ as follows. For any non-negative integer $n$ consider its finite $m$-ary expansion 
\[
n=\sum_{j=0}^{k_n} c_j(n)m^j;
\]
then we define
\[
t_n(m):= \sum_{j=0}^{k_n}c_j(n) \mod m.
\]

We refer to $\left(t_{n}(m)\right)_{n=0}^{\infty}$ as the \textbf{$\boldsymbol{TM_m}$ sequence}. Let $\Sigma\subset \Na$, $|\Sigma|=m$ and $f:\Sigma_m \to \Sigma$ a bijection, then $\left(f(t_n(m))\right)_{n=0}^{\infty}$ is called the $\boldsymbol{TM_m}$ \textbf{sequence over} $\mathbf{\Sigma}$.
\label{TMm01}
\end{def01}
\begin{obs01}
When $m=2$ we recover the classical Thue-Morse sequence (cf. \cite{queff}).
\end{obs01}
When the integer $m$ is understood, we will write $(t_n)_{n=0}^{\infty}$ in place of $(t_{n}(m))_{n=0}^{\infty}$.
\begin{propo01}
Fix an integer $m\geq 2$ and let $(t_n)_{n=0}^{\infty}$ be the associated $TM_m$ sequence. Then 
\begin{description}
\item[1. ] $t_n=t_{nm}$ for any $n\in \Na$, 
\item[2. ] $t_{n+1}-t_n\not\equiv 1 \mod m$ implies $n\equiv m-1 \mod m$,
\item[3. ] for any $r\in \Sigma_m$ and any $n \in \Na_{0}$, $t_{nm+r}\equiv t_{nm}+r \mod m$.

\end{description}
\label{obser}
\end{propo01}
\begin{proof}
If $n=[c_0\, c_{1}\, \ldots c_r]_m$, we have 
\[
n=\sum_{i=0}^{r} c_i m^i\; \Rightarrow nm=\sum_{i=0}^r c_i m^{i+1}\;  \Rightarrow \; nm=[0,\, c_0,\,c_{1},\,\ldots,\, c_r]_m.
\]
Suppose that $n\not\equiv m-1 \mod m$, then $n=[c_0\,c_{1}\,\ldots \,c_r]_m$ with $c_0\in \left\{0,1,\ldots,m-2\right\}$, therefore
\[
n+1=\sum_{i=0}^r c_im^i+1=(c_0+1)+c_1m+c_2m^2+\ldots+c_rm^r \; \Rightarrow \; n+1=[c_0+1,\,c_1,\,\ldots, \, c_r]_m.
\]
Hence,
\[
t_{n+1}\equiv (c_0+1)+c_1+\ldots+c_r \equiv 1+\sum_{i=0}^r c_i \equiv 1+t_n \mod m.
\]
The third part follows from the second, which states that $n\not\equiv m-1 \mod m$ implies $t_{n+1}-t_n\equiv 1 \mod m$. Since $nm\equiv 0 \mod m$, we have $t_{nm+1}-t_{nm}\equiv 1 \mod m$, $\ldots$, $t_{nm+(m-1)}-t_{nm + (m-2)}\equiv 1 \mod m$. By adding appropriately these equations we arrive at the desired conclusion.
\end{proof}
An immediate consequence of Proposition \ref{obser} is that $t_j=t_{j+1}=t_{j+2}$ is impossible. If this were the case, then $j\equiv -1 \mod m$ and $j+1 \equiv -1 \mod m$, which would imply $1\equiv 0 \mod m$.
\begin{def01}
Let $(x_i)_{i=0}^{\infty}$ be a sequence. If there exists an increasing sequence of positive numbers, $(n_j)_{j=1}^{\infty}$ such that for any positive integer $j$ and any $k\in \left\{0,\ldots,n_j\right\}$ the equality $x_k=x_{n_j-k}$ is satisfied, then $(x_i)_{i=0}^{\infty}$ is called \textbf{palindromic sequence}. If there exists an $N\in \Na$ such that $(x_i)_{i=N}^{\infty}$ is a palindromic sequence, then $(x_i)_{i=0}^{\infty}$ is called \textbf{eventually palindromic}.
\end{def01}

\begin{lem01}
Theorem \ref{Adabug} is still valid if \emph{palindromic} is replaced by \emph{eventually palindromic}. 
\end{lem01}
\begin{proof}
If $(a_{n})_{n=1}^{\infty}$ is a sequence of positive integers, $\alpha_n=[a_n;a_{n+1},a_{n+2},\ldots]$ for $n\geq 1$, and $\alpha=[0;a_1,a_2,a_3,\ldots]$, then for any $n\in \Na$ 
\[
T_{n-1}(\alpha_n)=\alpha, \; \text{ where } \; T_{n-1}(x)=\frac{p_{n-1}x-p_{n-2}}{q_{n-1}x-q_{n-2}}, 
\] 
and $(p_k/q_k)_{k=1}^{\infty}$ is the sequence of convergents of $\alpha$ (see \cite{khin}, p.8, eq.16 or \cite{schm}, p.18, proof of theorem 5C). The transformation $S_{n-1}=T^{-1}_{n-1}$ has integer coefficients because $|p_{n-2}q_{n-1}-p_{n-1}q_{n-2}|=1$ (see \cite{khin}, p. 5, Theorem 2), so $\alpha \in \QU(\alpha_n)$ and $\alpha_n\in \QU(\alpha)$ for all $n$. Hence, $\QU(\alpha)=\QU(\alpha_n)$ and $\alpha$ transcendental if and only if $\alpha_n$ is transcendental for any $n$. 
\end{proof}

Let $f:\Sigma_m \to \Na$ be injective. If $(x_n)_{n=0}^\infty$ is a sequence in $\Sigma_m$, then obviously $(f(x_{n}))_{n=0}^{\infty}$ is eventually periodic or eventually palindromic if and only if the original sequence is eventually periodic or eventually palindromic. The real numbers considered in Theorem \ref{Mien} have precisely the form $\alpha=[0;f(x_0),f(x_1),f(x_2),\ldots]$ with $x_n=t_{n}(m)$ and $f$ an injection as above. Therefore, it will be very important for us to verify that the $TM_m$ sequences are not eventually periodic. If this was the case, by Lagrange's theorem, all the numbers generated in this way would be quadratic irrationals.

\begin{teo01}
For any integer $m\geq 2$ the $TM_m$ sequence is not eventually periodic.
\label{aperiodic}
\end{teo01}
\begin{proof}

Let $m\in \Na$, $m\geq 2$, and let $(t_n)_{n=0}^{\infty}$ be the $TM_m$ sequence. It will be shown that there are not positive integers $a,b$ such that $t_{a+n}=t_{a+n+b}$ for all $n\in \Na$. Suppose otherwise. We may assume that $b$ is minimal. 
\begin{description}
\item \emph{Case $b<m$.} Suppose that $t_{a+n}=t_{a+b+n}$ for all $n\in \Na_{0}$ and some $a\in \Na_{0}$. In particular, if $a+a'\equiv 0 \mod m$, then Proposition \ref{obser}.3 yields
\[
t_{a+a'} \equiv t_{a+a'+b} \equiv t_{a+a'}+b \mod m
\]
which implies $b=0$, contradiction.
\item \emph{Case $b\equiv 0 \mod m$ and $b\geq m$.} The assumption implies $b=mb'$ for some positive integer $b'$. Then, by Proposition \ref{obser}.1., for any $n\in \Na$ we have
\[
t_{a+n}=t_{m(a+n)}=t_{m(a+n)+b}=t_{m(a+n+b')}=t_{a+n+b'},
\]
which contradicts the minimality of $b$.

\item \emph{Case $b\not\equiv 0 \mod m$ and $b>m$.} The closed interval $[a,a+2b]$ contains at least $2m+3$ consecutive integers, then there are $2$ multiples of $m$ strictly larger than $a$ in it, say $j=mq$ and $j+m=m(q+1)$. Suppose that $t_{j+m}-t_{j+m-1} \equiv 1 \mod m$ and $t_{j}-t_{j-1}\equiv 1 \mod m$, then $t_{j-m}=t_j=t_{j+m}$: indeed, by  Proposition \ref{obser},
\[
t_{(j-m)+r}=t_{j-m}+r, \; t_{j+r}=t_j+r, \; t_{j+m+r}=t_{j+m}+r  \quad \forall r\in \Sigma_m,
\]
so
\begin{align*}
t_j &\equiv 1+t_{j-1} \equiv 1+t_{(j-m)+(m-1)} \equiv 1+t_{j-m}+m-1 \equiv t_{j-m} \mod m, \nonumber\\
t_{j+m} &\equiv 1+t_{j+m-1} \equiv 1+t_{j+(m-1)} \equiv 1+t_{j}+m-1 \equiv t_{j} \mod m. \nonumber
\end{align*}
Then $t_{q-1}=t_q=t_{q+1}$, which is impossible by Proposition \ref{obser}.2. As a consequence, at least one of the expressions $t_{j+m}-t_{j+m-1}\not\equiv 1 \mod m$, $t_{j}-t_{j-1}\not\equiv 1 \mod m$ holds. We can assume, without losing generality, that $t_{j}-t_{j-1}\not\equiv 1 \mod m$ is true. Then $t_{j+b}-t_{j-1+b}\equiv t_{j}-t_{j-1}\not \equiv 1 \mod m$ which gives $j+b-1\equiv -1 \mod m$ and $j-1\equiv -1 \mod m$, implying that $b\equiv 0 \mod m$, contrary to hypothesis.
\end{description}
\end{proof}

If $\Sigma$ is a non-empty set, let $\Sigma^*$ be the set of all the finite words over $\Sigma$ and let $\Sigma^{\infty}\supset \Sigma^*$ be the set whose elements are the right infinite words over $\Sigma$ and the finite words over $\Sigma$. Morphisms constitute an important type of function mapping $\Sigma^*$ into $\Sigma^*$. The set $\Sigma^{*}$, including the empty word, equipped with the operation of concatenation, forms a monoid.

\begin{def01}
A \textbf{morphism} is a monoid homomorphism, i.e.
\[
\bfw=w_0w_1\ldots w_k\in \Sigma^*  \; \Rightarrow \; \vphi(\bfw)=\vphi(w_0)\vphi(w_1)\ldots \vphi(w_k).
\]
Any morphism $\vphi$ can be extended to a function $\overline{\vphi}:\Sigma^{\infty} \to \Sigma^{\infty}$ by
\[
\overline{\vphi}(\bfw)=\begin{cases}
\vphi(w_0)\vphi(w_1)\vphi(w_2)\ldots \; \text{ if } \; \bfw=w_0 w_1 w_2 \ldots \in \Sigma^{\infty}\setminus \Sigma^*, \\
\vphi(\bfw) \; \text{ if } \; \bfw\in \Sigma^*.
\end{cases} 
\]
In this case, $\overline{\vphi}$ is said to be \textbf{generated} by $\vphi$.
\label{morph01}
\end{def01}

Clearly, any morphism $\vphi:\Sigma^* \to \Sigma^*$, and hence the morphism generated by it, is completely determined by $\vphi(j)$ with $j$ ranging over $\Sigma$. In what follows we will use the same notation for $\vphi$ and $\overline{\vphi}$, also $\vphi^1=\vphi$ and for any $n\in \Na$ we will define $\vphi^{n+1}:=\vphi^n\circ \vphi$.

\begin{def01}
A morphism $\vphi:\Sigma^*\to \Sigma^*$ is \textbf{prolongable on} $\boldsymbol{j}\in \Sigma$ if $\vphi(j)=j\bfx$ for some $\bfx\in\Sigma^*$. If $\vphi$ is prolongable on $j$, then the \textbf{fixed point of $\boldsymbol{\overline{\vphi}}$ based on $\boldsymbol{j}$} is the word $j\bfx\vphi(\bfx)\vphi^2(\bfx)\vphi^3(\bfx)\ldots$. 
\label{morph02}
\end{def01}
\begin{def01}
Let $\vphi:\Sigma^*\to \Sigma^*$ be a morphism and $k\in \Na$. Then $\vphi$ is a $\boldsymbol{k}$\textbf{-uniform} morphism if $|\vphi(j)|=k$ for any $j\in \Sigma$.
\end{def01}

\begin{teo01}
Let $m\geq 2$ be an integer and let $(t_n)_{n=0}^{\infty}$ be its associated $TM_m$ sequence. If $m=2$, then $(t_{n})_{n=0}^{\infty}$ is palindromic. If $m>2$, then $(t_{n})_{n=0}^{\infty}$ is not eventually palindromic.
\label{nopalin}
\end{teo01}
\begin{proof}
For the case $m=2$ we give the argument found in \cite{allouche} as it will motivate ideas which will be used in the sequel. Applied to the Thue-Morse sequence ($m=2$), Proposition \ref{obser} implies that $(t_n)_{n=0}^{\infty}$ is palindromic, where the sequence $(n_j)_{j=1}^{\infty}$ is given by $n_j=4^{j-1}$. To see this, let $\vphi:\left\{0,1\right\}^{\infty}\to \left\{0,1\right\}^{\infty}$ be the morphism generated by $\vphi(0)=01$ and $\vphi(1)=10$. By definition, $\vphi^2(0)=0110$, $\vphi^2(1)=1001$, thus the set of palindromic words over $\left\{0,1\right\}$ is invariant under $\vphi^2$. Indeed, let $\bfw\in\Sigma^*$ be a palindrome, $\bfw=w_0w_1\ldots w_{k-1}$, then, if $\overline{w}_j:=1-w_j$, we have
\[
\vphi^2(\bfw)=\vphi^2(w_0)\vphi^2(w_1)\ldots\vphi^2(w_{k-1})=w_0\,\overline{w}_0\,\overline{w}_0w_0\, w_1\,\overline{w}_1\,\overline{w}_1\,w_1\ldots w_{k-1}\, \overline{w}_{k-1}\,\overline{w}_{k-1}\,w_{k-1}. 
\]
Noting that $|\vphi(\bfw)|=4k$, let $j\in \left\{0,1,\ldots,4k-1\right\}$ and let $q,r\in\Na_0$ be such that $j=4q+r$ and $0\leq r < 4$. Clearly,
\[
\vphi(\bfw)_j=\vphi(\bfw)_{4q+r}=\begin{cases} w_q \text{ if } r\in \left\{0,3\right\}, \\ \overline{w}_q \text{ if } r\in \left\{1,2\right\}.\end{cases}
\]
Therefore, since $r\in \left\{0,3\right\}$ if and only if $3-r\in \left\{0,3\right\}$ and $r\in \left\{1,2\right\}$ if and only if $3-r\in \left\{1,2\right\}$,
\[
\vphi(\bfw)_{4k-1-j}=\vphi(\bfw)_{4k-1-4q-r}=\vphi(\bfw)_{4(k-1-q)+(3-r)}=\begin{cases} w_{k-1-q} \text{ if } r\in \left\{0,3\right\},\\ \overline{w}_{k-1-q} \text{ if } r\in \left\{1,2\right\}. \end{cases}
\]
So, if $\bfw$ is a palindrome, for any $l\in \left\{0,1,\ldots,k-1\right\}$ the equality $w_l=w_{k-1-l}$ holds implying that $\vphi(\bfw)_j=\vphi(\bfw)_{4k-1-j}$ holds too and, consequently, that $\vphi(\bfw)$ is a palindrome. By definition, $\vphi^2(0)=0110$ is a palindrome and by induction $\vphi^{2n}(0)$ is a palindrome for any $n\in \Na$. Note that $\vphi$ is prolongable on $0$. Moreover, we have $\vphi(0)=01$, $\vphi^2(0)=0110=01\vphi(1)$, if $\vphi^{2n}(0)=01\vphi(1)\ldots \vphi^{2n-1}(0)$, then
\begin{align*}
\vphi^{2(n+1)}(0)&=\vphi^2(\vphi^n(0))=\vphi^2(01\vphi(1)\ldots\vphi^{2n-1}(0))=\vphi^2(0)\vphi^2(1)\ldots\vphi^{2n+1}(0) \nonumber\\
&=01\vphi(1)\ldots\vphi^{2n+1}(0). \nonumber\\
\end{align*}
This means that $\vphi^{2n}(0)$ is an initial segment of the fixed point of $\vphi$ based on $0$, $\bfv$. In particular, $\bfv$ is palindromic where the initial palindromic segments are given by $\vphi^{2n}(0)$. By Theorem 3.1 of the next section \ref{defieq} (whose proof does not depend on the current Theorem), $\bfv$ is precisely $TM_2$. Therefore, $TM_2$ is a palindromic sequence.

Now we consider the case $m\geq 3$. In $(t_n)_{n=0}^{\infty}$ the word $011$ appears infinitely often, for
\[
j:=\big[m-2,\underbrace{m-1,\ldots,m-1}_\textrm{m-2 \text{ times}}\big]_m, j+1=\big[\underbrace{m-1,\ldots,m-1,m-1}_\textrm{m-1 \text{ times}}\big]_m, j+2=\big[\underbrace{0,\ldots,0}_\textrm{m-1 \text{ times}},1\big]_m
\]
give $t_j=0$, $t_{j+1}=1$, $t_{j+2}=1$. We observe that $t_{j'}=0$, $t_{j'+1}=1$, $t_{j'+2}=1$ when $j_m'$ is defined by
\[
j':=\big[m-2,\underbrace{m-1,\ldots,m-1}_\textrm{m-2 \text{ times}},0,m-2,2\big]_m=j+(m-2)m^m+2m^{m+1}.
\]

In general, $j(k):=j+(m-2)m^k+2m^{k+1}$ for $k\geq m$ give rise to an occurence of $011$; that is,  $t_{j(k)}=0$, $t_{j(k)+1}=1$, $t_{j(k)+2}=1$ where
\[
j(k):=\big[m-2,\underbrace{m-1,\ldots,m-1}_\textrm{m-2 \text{ times}},\underbrace{0,\ldots,0}_\textrm{k+1-m \text{ times}},m-2,2\big]_m=j+(m-2)m^k+2m^{k+1}.
\]

However, the configuration $110$ is impossible. For suppose we had $t_k=1$, $t_{k+1}=1$, $t_{k+2}=0$ for some $k\in \Na_0$. Since $m\geq 3$, $-1\not\equiv 1 \mod m$ and so, $t_{k+2}-t_{k+1}\not\equiv 1 \mod m$. We also have $t_{k+1}-t_k\not \equiv 1 \mod m$. By Proposition \ref{obser}, $k+1\equiv m-1 \mod m$ and $k\equiv m-1 \mod m$. The last congruences yield $1\equiv 0 \mod m$, which is absurd. Hence, the $TM_m$ sequence is not eventually palindromic.
\end{proof}

\section{An alternative definition of the $TM_m$ sequences.}
The classical Thue-Morse sequence ($TM_2$) may be defined in several ways. It can be defined as the fixed point based on $0$ of a certain morphism (see \cite{allouche}, Corollary 1.7.7, p.23). In this spirit, for any fixed $m\in \Na$, $m\geq 2$, we define the morphism $\vphi:\Sigma_m^*\to \Sigma_m^*$ by 
\[
\vphi(j):=j, \, j+1,  \, \ldots,  \, j+m-1  \quad \forall j \in \Sigma_m.
\]
Addition in the above definition must be understood modulo $m$. For instance, if $m=4$ and $j=2$, then $\vphi(j)=2 \, 3  \, 0  \, 1$ (as noted earlier, $\vphi$ is completely determined by its action on $\Sigma_m$). By definition, $\vphi$ is prolongable on $0$ and, in fact, on any element of $\Sigma_m$. Also, $\vphi^k$ is $m^k$-uniform and for any $k\in \Na$.

\begin{def01}
For any integer $m\geq 2$ let $\vphi$ be the morphism defined as above. The $\boldsymbol{TM_m}$ \textbf{sequence} is the fixed point of $\vphi$ based at $0$.
\label{TMm02}
\end{def01}
From now on, we shall asume that $m$ has been fixed.
\begin{teo01}
Definitions \ref{TMm01} and \ref{TMm02} are equivalent.
\label{defieq}
\end{teo01}
\begin{lem01}
For any $c_0,\ldots,c_n,c_{n+1} \in \Sigma_m$ define the word $\bfc:=c_0\ldots c_{n+1}$. Then, the following equality holds
\begin{equation}
\vphi^{n+1}(c_{n+1})_{\left[\bfc\right]_m}=\sum_{i=0}^{n+1} c_i \mod m.
\label{fimorf}
\end{equation}
\label{lemcasi}
\end{lem01}
We will need to deal with long expressions for words, so it will be convenient to adopt new notation. If $\bfw\in \Sigma_m^*$, then $\left\{ \bfw\right\}_j:=w_j$ for any $j\in \left\{0,1,\ldots,|\bfw|-1\right\}$.
\begin{proof}
Consider first the case $n=0$. Let $c_0$, $c_1 \in \Sigma_m$, then 
\[
\vphi(c_1)_{c_0}=\left\{ c_1,\, c_1+1,\, \ldots\,,c_1+m-1\right\}_{c_0}=c_1+c_0.
\]
Now consider the case $n=1$: since $\vphi$ is an $m$-uniform morphism and by definition of $\vphi$,
\[
\vphi^2(c_2)_{c_1m+c_0}=\left\{\vphi(c_2)\vphi(c_2+1)\ldots\vphi(c_2+m-1)\right\}_{c_1m+c_0}=\vphi(c_2+c_1)_{c_0}=c_2+c_1+c_0.
\]
Suppose that \eqref{fimorf} holds for $n-1$. Taking $c_0,c_1,\ldots,c_{n+1} \in \Sigma_m$, and $\bfc:=c_0c_1\ldots c_n$,
\begin{align*}
\vphi^{n+1}(c_{n+1})_{[\bfc]_m}&=\left\{\vphi^n(c_{n+1}),\vphi^n(c_{n+1}+1),\ldots, \vphi^n(c_{n+1}+m-1)\right\}_{[\bfc]_m} \nonumber\\
&=\vphi^n(c_{n+1}+c_n)_{\sum_{i=0}^{n-1} c_i m^i} \nonumber\\
&=\sum_{i=0}^{n+1} c_i. \nonumber\\
\end{align*}
The second equality follows from $\vphi^n$ being $m^n$-uniform morphism, and the third equality is valid by induction. 
\end{proof}
\begin{proof}[Proof of Theorem \ref{defieq}]
Let $\bfv \in \Sigma_m^{\infty}$ be the fixed point of $\vphi$ based at $0$, then
\begin{align*}
\bfv &= 0,\,1,\,\ldots, \,m-1, \vphi(1,\,\ldots, \,m-1),\vphi^2(1,\,\ldots, \,m-1), \ldots \nonumber\\
&=\vphi(0)\,\vphi(1)\,\ldots\, \vphi(m-1) \vphi^2(1,\,\ldots,\, m-1)\vphi^3(1,\, \ldots, \,m-1)\ldots \nonumber\\
&=\vphi^2(0)\vphi^2(1)\ldots \vphi^2(m-1)\vphi^3(1,\,\ldots,\,m-1)\vphi^4(1,\,\ldots,\,m-1)\ldots \nonumber\\
&\vdots \nonumber\\
&=\vphi^n(0)\vphi^n(1)\ldots \vphi^n(m-1)\vphi^{n+1}(1,\,\ldots,\,m-1)\vphi^{n+2}(1,\,\ldots,\,m-1)\ldots. \nonumber\\
\end{align*}
Let $j$ be nonnegative, $j=\sum_{i=0}^n c_i m^i$, and define $\bfc:=(j)_m$. Taking a sufficiently long initial segment and applying Lemma \ref{lemcasi} we get
\[
v_j=\vphi^{n+1}(0)_{j}=\vphi^{n+1}(0)_{[\bfc]_m}=\sum_{i=0}^n c_i \mod m=t_j.
\]
\end{proof}

\section{Proof of Theorem \ref{Mien}. }
Definition \ref{TMm02} makes available to us the theory of word morphisms. In particular, certain extremely useful combinatorial aspects will be employed. 
\begin{def01}
Let $\mathbf{w}=a_0a_1a_2\ldots$ be an infinite word over a finite alphabet $\Sigma$. The \textbf{complexity function}, $p_{\mathbf{w}}:\Na\to \Na$, is given by 
\[
p_{\mathbf{w}}(n):=\big| \left\{ \mathbf{y} \in \Sigma^*: (\mathbf{y} \text{ is a subword of }\mathbf{w}) \, \& \, (|\mathbf{y}|=n) \right\}\big|\leq |\Sigma|^n.
\]
\end{def01}

\begin{lem01}
Let $n,m$ be positive integers with $m\geq 2$ and let $(t_j)_{j=1}^{\infty}$ the $TM_m$ sequence. If $\bft=t_{0} t_{1} t_{2} \ldots$, then $p_{\mathbf{t}}(n)\leq m^3 n$. 
\label{lemacasi}
\end{lem01}
\begin{proof}
Let $n$ be a fixed positive integer and $r\in \Na$ such that $m^{r-1}\leq n <m^r$, and let $\bft[k]$ be the set of subwords of $\bft$ of length $k$. Let $f:\Sigma_{m^r}\times \bft[2]\to \bft[n]$ be a function given by
\begin{equation}
f(s,\bfv):=\vphi^r(\bfv)_{[s:s+n]}.
\label{defdeefe}
\end{equation}
Since $\bft$ is a fixed point of $\vphi$ and $\vphi$ is an $m$-uniform morphism, $\bft$ is a fixed point of $\vphi^p$ and $\vphi^p$ is an $m^p$-uniform morphism for any $p\in \Na$. Therefore, we have
\[
\bft=t_0t_1t_2\ldots=\vphi^p(t_0)\vphi^p(t_1)\vphi^p(t_2)\ldots=\vphi^p(\bft) \quad \forall \, p\in \Na.
\]
Let $t_jt_{j+1}$ be an occurrence of $\bfv$. The equation above and the just stated properties of $\vphi$ imply
\[
\vphi^r(\bfv) = \vphi^r(t_j)\vphi^r(t_{j+1})=t_{jm^r}t_{jm^r+1}\ldots t_{(j+1)m^r-1}t_{(j+1)m^r}t_{(j+1)m^r+1}\ldots t_{(j+2)m^r-1}, 
\]
and we obtain
\[
\vphi^r(\bfv)_{[s:s+n]}=t_{jm^r+s}t_{jm^r+s+1}\ldots t_{jm^r+s+n-1}.
\]

Let $\bfx$ be any subword of $\bft$ of length $n$, then $\bfx=t_i\ldots t_{i+n-1}$ for some $i$. By the Euclidean Algorithm, there are two nonnegative integers $q',s'$ such that $i=q'm^r+s'$ and $0 \leq s'<m^r$. Obviously, $\bfx$ is a subword of 
\[
\vphi^r(t_{q'}t_{q'+1})=t_{q'm^r}\ldots t_{(q'+1)m^r-1}t_{(q'+1)m^r}\ldots t_{(q'+2)m^r-1}; 
\]
in fact, $\bfx=f(s',\bfv)$ with $\bfv=t_{q'}t_{q'+1}$ which proves that $f$ is onto. Therefore 
\[
p_{\bft}(n)=\big|\bft[n]\big|\leq \big|\Sigma_{m^r}\times \bft[2]\big|=m^rm^2=m^{r-1}m^3\leq nm^3
\]
which yields
\begin{equation}
\limsup_{n\to +\infty} \frac{p_{\bft}(n)}{n} \leq m^3.
\label{lisumTM}
\end{equation}
\end{proof}

\begin{obs01}
Our argument is a special case of the proof of a more general theorem (see \cite{allouche}, Theorem 10.3.1, p. 304). The next theorem, whose proof is ommitted, along with Lemma \ref{lemacasi}, will lead us to Theorem \ref{Mien}. 
\end{obs01}

\begin{teo01}[Bugeaud] 
Let $\mathbf{a}=(a_n)_{n=1}^{\infty}$ be a sequence of positive integers which is not eventually periodic. If $\alpha:=[0,a_1,a_2,a_3,\ldots]$ is algebraic and $p_n(\mathbf{a})$ is the complexity function, then
\begin{equation}
\lim_{n\to +\infty} \frac{p_n(\mathbf{a})}{n}=+\infty.
\end{equation}
\label{Bug}
\end{teo01}
\begin{proof}
See \cite{bugo02}, Theorem 1.1.
\end{proof}

\begin{flushleft}
\emph{Proof of Theorem \ref{Mien}}.
\end{flushleft}

Let $m\geq 3$ be an integer, $(t_j)_{j=0}^{\infty}$ the $TM_m$ sequence, $\Sigma\subset \Na$, $\psi:\Sigma_m \to \Sigma$ bijective, $\bfw:=\psi(t_{0})\psi(t_{1})\psi(t_{2})\ldots$ and $\alpha=[0;\psi(t_{0}),\psi(t_{1}),\psi(t_{2}),\ldots]$. The conclusion of Lemma \ref{lemacasi} remains valid--$\psi$ is a bijection-- when applying it to $\bfw$, so \eqref{lisumTM} is valid and, by Theorem \ref{aperiodic}, $\bfw$ is not eventually periodic. Since 
\[
\limsup_{n\to+\infty} \frac{p_n(\bfw)}{n}<\infty,
\]
the real number $\alpha$ is transcendental.

\end{document}